\documentclass[12pt,oneside,english]{amsart}
\usepackage[T1]{fontenc}
\usepackage[latin9]{inputenc}
\usepackage{babel}
\usepackage{amstext}
\usepackage{amsthm}
\usepackage{amssymb}
\usepackage[unicode=true,pdfusetitle,
 bookmarks=true,bookmarksnumbered=false,bookmarksopen=false,
 breaklinks=false,pdfborder={0 0 1},backref=false,colorlinks=false]
 {hyperref}

\makeatletter
\numberwithin{equation}{section}
\theoremstyle{plain}
\newtheorem{conjecture}{\protect\conjecturename}[section]
\theoremstyle{plain}
\newtheorem{thm}{\protect\theoremname}[section]
\theoremstyle{plain}
\newtheorem{lem}{\protect\lemmaname}[section]

\usepackage{textcomp}

\usepackage{listings}
\usepackage{xcolor}

\lstdefinelanguage{Maple}{
  morekeywords={proc,local,for,to,do,end,while,if,then,elif,else,error,
    printf,print,restart,seq,add,expand,diff,degree,coeff,eval,subs,min,max},
  sensitive=true,
  morecomment=[l]{\#},
  morestring=[b]",
  alsoletter={_}
}

\makeatother

\providecommand{\conjecturename}{Conjecture}
\providecommand{\lemmaname}{Lemma}
\providecommand{\theoremname}{Theorem}

\begin{document}
\title[a conjecture of Ruscheweyh, Salinas and Sugawa]{Disproof of a conjecture of Ruscheweyh et al. on universal convexity
of Ramanujan generating function}
\author{Bing He}
\address{School of Mathematics and Statistics, HNP-LAMA, Central South University
\\
 Changsha 410083, Hunan, People's Republic of China}
\email{yuhelingyun@foxmail.com; yuhe123456@foxmail.com}
\begin{abstract}
We disprove a conjecture of Ruscheweyh, Salinas and Sugawa concerning
the universal convexity of the normalized generating function associated
with the Ramanujan sequence $\theta(n).$ More precisely, we prove
that the function
\[
\sigma(z)=\frac{1}{\theta(1)-1/3}\sum_{n\geq1}(\theta(n)-1/3)z^{n}
\]
is not universally convex on the slit domain $\mathbb{C}\backslash[1,\infty).$
The proof uses Watson\textquoteright s integral representation of
$\theta(n),$ asymptotic estimates for an auxiliary function $h,$
and the representation of derivatives of universally convex functions
via logarithmic integrals. As a general extension of this result,
we also show that the generalized Ramanujan generating function 
\[
\sigma_{\alpha}(z)=\frac{1}{\theta(1)-1/3}\sum_{n\geq1}(\theta(n)-1/3)\frac{z^{n}}{n^{\alpha}}
\]
is universally convex if and only if $\alpha\geq1.$
\end{abstract}

\maketitle

\section{Introduction}

For $n=0,1,2,\ldots$ the Ramanujan sequence $\{\theta(n)\}_{n\geq0}$
is defined by 
\[
\frac{e^{n}}{2}=\sum_{k=0}^{n-1}\frac{n^{k}}{k!}+\frac{n^{n}}{n!}\theta(n)\quad(n\geq1),\qquad\theta(0)=\frac{1}{2}.
\]
A famous problem of Ramanujan asks whether 
\begin{equation}
\frac{1}{3}\leq\theta(n)\leq\frac{1}{2}\qquad(n\geq0).\label{eq:1-1}
\end{equation}
This was refined by Ramanujan in his first letter to Hardy \cite{R}
as 
\begin{equation}
\theta(n)=\frac{1}{3}+\frac{4}{135}\cdot\frac{1}{n+k_{n}},\quad k_{n}\in\left[\frac{2}{21},\frac{8}{45}\right],\quad n\geq0.\label{eq:1-2}
\end{equation}
The first proofs of \eqref{eq:1-1} were given by Szegö \cite{S}
and Watson \cite{W} while a proof of \eqref{eq:1-2} was provided
by Flajolet et al. \cite{FGKP}. Adell and Jodrá \cite{AJ} showed
that there exists a probability distribution function $F$ on {[}0,1{]}
such that 
\[
\theta(n)-\frac{1}{3}=\frac{1}{6}\int_{0}^{1}x^{n}dF(x),\qquad n\geq0.
\]
This, together with the Hausdorff moment theorem, implies that the
sequence $\{\theta(n)\}_{n\geq0}$ is completely monotone\footnote{A sequence $\{\alpha(n)\}$ is said to be completely monotone if for
all $r,n,$
\[
\sum_{j=0}^{r}(-1)^{j}\binom{r}{j}\alpha(n+j)\geq0.
\]
}. In \cite{K} Koumandos used 
\[
\theta(n)=1+\frac{n}{2}\left[\int_{0}^{1}\left(ue^{1-u}\right)^{n}du-\int_{1}^{\infty}\left(ve^{1-v}\right)^{n}dv\right]
\]
as the definition of the function $\theta(n)$ for all positive real
numbers $n$ and proved that the function 
\[
\frac{4}{135}-x\left(\theta(x)-\frac{1}{3}\right)
\]
 is completely monotonic on $[0,\infty).$ Consequently, the sequence
\[
\left\{ \frac{4}{135}-n\left(\theta(n)-\frac{1}{3}\right)\right\} _{n\geq0}
\]
is completely monotone. In 2018, Bakan, Ruscheweyh and Salinas \cite{BRS}
derived that there exists a continuous probability distribution function
$G$ on {[}0,1{]} such that 
\[
(n+1)\left(\theta(n)-\frac{1}{3}\right)=\frac{4}{135}+\frac{37}{270}\int_{0}^{1}x^{n}dG(x),\quad n\geq0,
\]
which implies that the sequence 
\[
\left\{ (n+1)\left(\theta(n)-\frac{1}{3}\right)\right\} _{n\geq0}
\]
is completely monotone.

Let $\Lambda$ denote the slit domain $\mathbb{C}\backslash[1,\infty).$
A function $f$ analytic on $\Lambda$ is called normalized if 
\[
f(0)=0,\qquad f'(0)=1.
\]
A normalized analytic function $f$ on $\ensuremath{\Lambda}$ is
said to be universally convex if it maps every circular domain contained
in $\Lambda$ conformally onto a convex domain. Here, by a circular
domain we mean an open disk or an open half-plane in $\mathbb{C}.$

In \cite[p. 294]{RSS}, Ruscheweyh, Salinas and Sugawa proposed the
following conjecture.
\begin{conjecture}
The function
\[
\sigma(z):=\frac{1}{\theta(1)-1/3}\sum_{n\geq1}(\theta(n)-1/3)z^{n}
\]
is universally convex.
\end{conjecture}
In this paper, we disprove this conjecture by proving the following
result.
\begin{thm}
\label{t1} The function $\sigma(z)$ is not universally convex.
\end{thm}
The remainder of the paper is organized as follows. Section \ref{sec:Au}
proves Theorem \ref{t1} by combining Watson\textquoteright s integral
representation of $\theta(n),$ asymptotic estimates for an auxiliary
function $h$ and the representation of derivatives of universally
convex functions via logarithmic integrals. Section \ref{sec:A-p}
establishes a positivity lemma for $h,$ which is the key technical
ingredient in Section \ref{sec:Univ}. In Section \ref{sec:Univ},
as a general extension of Theorem \ref{t1}, we study the generalized
Ramanujan generating function $\sigma_{\alpha}(z)$ and show that
it is universally convex if and only if $\alpha\geq1.$

\section{\label{sec:Au} Proof of Theorem \ref{t1}}

\subsection{An auxiliary result}

For a probability measure $\mu$ on {[}0, 1{]}, set

\[
L_{\mu}(z):=\int_{[0,1]}\log\frac{1}{1-tz}d\mu(t),\quad z\in\Lambda
\]
where the logarithm is the analytic branch equal to 0 at $\ensuremath{z=0}.$

 A normalized analytic function $f$ on $\ensuremath{\Lambda}$ is
universally starlike if it maps every circular domain $\Omega\subset\Lambda$
containing 0 conformally onto a domain starlike with respect to 0.
Let $\mathcal{R}_{\alpha}^{u}$ denote the set of universally prestar-like
functions of order $\alpha.$
\begin{lem}
\label{l1} A normalized analytic function $f$ on $\Lambda$ is universally
convex if and only if 
\[
f'(z)=\exp\bigl(2L_{\mu}(z)\bigr)
\]
 for some probability measure $\mu$ on $\ensuremath{[0,1]}.$
\end{lem}
\begin{proof}
Suppose that the normalized analytic function $f$ is universally
convex. By \cite[Theorem 1.6]{RSS}, we know that the function $f$
is universally prestarlike of order $\ensuremath{0},$ i.e. $f\in\mathcal{R}_{0}^{u}.$
Then, by \cite[Corollary 1.1]{RSS}, there exists a probability measure
$\mu$ on $[0,1]$ such that 
\[
\frac{f(z)}{z}=I_{2}(z)*\exp\left(\int_{0}^{1}\log\frac{1}{(1-tz)^{2}}\,d\mu(t)\right),
\]
 where $*$ denotes the Hadamard product and 
\[
I_{\gamma}(z):={}_{2}F_{1}(1,1,\gamma,z).
\]
Since 
\[
\log\frac{1}{(1-tz)^{2}}=2\log\frac{1}{1-tz},
\]
 we obtain 
\[
\frac{f(z)}{z}=I_{2}(z)*\exp\bigl(2L_{\mu}(z)\bigr).
\]
 Let 
\[
H(z)=\exp\bigl(2L_{\mu}(z)\bigr)=\sum_{n=0}^{\infty}h_{n}z^{n}.
\]
Note that
\[
I_{2}(z)={}_{2}F_{1}(1,1,2,z)=\sum_{n=0}^{\infty}\frac{z^{n}}{n+1}=-\frac{\log(1-z)}{z}.
\]
 Then 
\[
(I_{2}*H)(z)=\sum_{n=0}^{\infty}\frac{h_{n}}{n+1}z^{n}.
\]
 Thus 
\[
f(z)=z(I_{2}*H)(z)=\sum_{n=0}^{\infty}\frac{h_{n}}{n+1}z^{n+1}=\sum_{m=1}^{\infty}\frac{h_{m-1}}{m}z^{m}.
\]
 Differentiating termwise gives
\[
f'(z)=\sum_{m=1}^{\infty}h_{m-1}z^{m-1}=\sum_{n=0}^{\infty}h_{n}z^{n}=H(z).
\]
 Therefore 
\[
f'(z)=H(z)=\exp\bigl(2L_{\mu}(z)\bigr).
\]

Conversely, suppose 
\[
f'(z)=\exp\bigl(2L_{\mu}(z)\bigr)=H(z).
\]
Since $\ensuremath{f(0)=0},$ we have
\[
f(z)=\int_{0}^{z}H(s)\,ds=\sum_{n=0}^{\infty}\frac{h_{n}}{n+1}z^{n+1}=z(I_{2}*H)(z).
\]
 Hence 
\[
\frac{f(z)}{z}=I_{2}(z)*\exp\bigl(2L_{\mu}(z)\bigr).
\]
 By \cite[Corollary 1.1]{RSS} we know that $\ensuremath{f\in\mathcal{R}_{0}^{u}},$
which together with \cite[Theorem 1.6]{RSS}, implies that $f$ is
universally convex. This proves Lemma \ref{l1}.
\end{proof}

\subsection{\label{subsec:PT1} Proof of Theorem \ref{t1}}

Before proving Theorem \ref{t1}, we introduce the auxiliary functions
$u(x),$ $U(x),$ and $h(x),$ which encode the two real branches
of the equation $ue^{1-u}=e^{-x}.$ These functions will allow us
to express $\sigma'(-r)$ as an integral involving a positive kernel
and a weight function $h.$ The estimates for $h$ obtained below
are the key to comparing the actual decay of $\sigma'(-r)$ with the
decay forced by the representation in Lemma \ref{l1}.

Let $u(x)$ and $U(x)$ be the two real solutions of
\begin{equation}
u(x)e^{1-u(x)}=U(x)e^{1-U(x)}=e^{-x},\quad0<u(x)<1<U(x)\quad(x>0).\label{eq:3-1}
\end{equation}
The function $U$ is strictly increasing on $(0,\infty)$ and satisfies
$U(x)\rightarrow\infty$ as $x\rightarrow\infty,$ whereas $u$ is
strictly decreasing on $(0,\infty)$ and satisfies $u(x)\rightarrow0$
as $x\rightarrow\infty.$

Set 
\[
h(x):=-\frac{1}{2}\left(U^{\prime\prime}(x)+u^{\prime\prime}(x)\right).
\]
Then, by \cite[Theorem 4 and eq.(24)]{BRS},
\[
h(x)>0,\quad h'(x)<0,\quad h(0+)=\frac{4}{135}.
\]

We now divive the proof of Theorem \ref{t1} into several steps.

Step 1: Representation of $\sigma^{\prime}(-r).$ It follows from
\eqref{eq:3-1} that 
\[
\ensuremath{U-\log U=x+1}.
\]
Differentiating both sides of this identity gives
\[
U'-\frac{U'}{U}=1,\qquad U'=\frac{U}{U-1}.
\]
Similarly, we have
\[
u'=\frac{u}{u-1}.
\]
Differentiating again gives
\begin{equation}
U^{\prime\prime}(x)=-\frac{U(x)}{(U(x)-1)^{3}},\quad u^{\prime\prime}(x)=\frac{u(x)}{(1-u(x))^{3}}.\label{eq:3-5}
\end{equation}

From \eqref{eq:3-1} we know that $\ensuremath{U(x)/x\rightarrow1},$
as $x\rightarrow\infty.$ Also, $\ensuremath{u(x)=O\left(e^{-x}\right)}.$
Thus \eqref{eq:3-5} gives
\begin{equation}
h(x)\sim\frac{1}{2x^{2}}\quad(x\rightarrow\infty).\label{eq:2-6}
\end{equation}

Set 
\[
a_{n}:=\theta(n)-\frac{1}{3},\quad c:=a_{1}=\frac{e}{2}-\frac{4}{3}>0.
\]
In \cite[p. 297]{W} Watson obtained that
\begin{equation}
a_{n}=\int_{0}^{\infty}e^{-nx}h(x)dx,\quad n\geq0.\label{eq:2-1}
\end{equation}
Then

\[
\int_{0}^{\infty}h(x)dx=a_{0}=\frac{1}{6}.
\]

Summing \eqref{eq:2-1} over all $n\geq0$ first for $\ensuremath{|z|<1},$
then continuing analytically to $\ensuremath{\Lambda},$ gives
\begin{equation}
\sigma(z)=\frac{1}{c}\int_{0}^{\infty}\frac{ze^{-x}}{1-ze^{-x}}h(x)dx.\label{eq:2-5}
\end{equation}
Differentiation under the integral is valid on compact subsets of
$\ensuremath{\Lambda}.$ In particular, for $\ensuremath{r>0},$
\[
\sigma^{\prime}(-r)=\frac{1}{c}\int_{0}^{\infty}\frac{e^{-x}h(x)}{\left(1+re^{-x}\right)^{2}}dx.
\]

Step 2: Upper bound for $Y(r).$ Let
\[
Y(r):=r\sigma^{\prime}(-r)=\frac{1}{c}\int_{0}^{\infty}\frac{re^{-x}h(x)}{\left(1+re^{-x}\right)^{2}}dx
\]
and 
\[
K_{r}(x):=\frac{re^{-x}}{\left(1+re^{-x}\right)^{2}}.
\]
The kernel $K_{r}(x)$ satisfies
\[
0\leq\frac{re^{-x}}{\left(1+re^{-x}\right)^{2}}\leq\frac{1}{4}\quad(r>0,x\geq0)
\]
and tends to 0 for each fixed $\ensuremath{x}$ as $r\rightarrow\infty.$
Then
\[
Y(r)\leq\frac{1}{4c}\int_{0}^{\infty}h(x)dx=\frac{1}{24c}.
\]
Since $\int_{0}^{\infty}h(x)dx=\frac{1}{6}<\infty,$ the dominated
convergence gives
\begin{equation}
Y(r)\rightarrow0\quad(r\rightarrow\infty).\label{eq:2-4}
\end{equation}

Step 3: Lower bound for $Y(r).$ For the lower bound, restrict the
integral to the interval $\ensuremath{x\in[\log r,\log r+1]}.$ On
this interval, the quantity $q:=re^{-x}$ lies in $\ensuremath{\left[e^{-1},1\right]}.$
Therefore 
\[
\frac{re^{-x}}{(1+re^{-x})^{2}}\ge\min_{q\in[e^{-1},1]}\frac{q}{(1+q)^{2}}=\frac{e}{(e+1)^{2}}=:m>0.
\]
Also, since $h(x)\sim\frac{1}{2x^{2}}$ for large $x,$ we have $h(x)\ge\frac{1}{4x^{2}}.$
On $\ensuremath{[\log r,\log r+1]},$ we have
\[
h(x)\ge\frac{1}{4(\log r+1)^{2}}.
\]
Hence, for sufficiently large $\ensuremath{r},$ 
\[
Y(r)\ge\frac{1}{c}\int_{\log r}^{\log r+1}mh(x)\,dx\ge\frac{m}{4c(\log r+1)^{2}}\ge\frac{m}{16c(\log r)^{2}}.
\]

Step 4: Logarithmic asymptotics of $\sigma^{\prime}(-r).$ Combining
the upper bound and the lower bound for $Y(r),$we obtain that, for
all sufficiently large $\ensuremath{r},$
\begin{equation}
\frac{C_{1}}{(\log r)^{2}}\le Y(r)\le C_{2}\label{eq:3-6}
\end{equation}
 with 
\[
C_{1}=\frac{e}{16c(e+1)^{2}}>0,\qquad C_{2}=\frac{1}{24c}>0.
\]
Note that $Y(r)=r\sigma^{\prime}(-r).$ Taking logarithms of \eqref{eq:3-6}
and dividing by $\log r$ gives
\begin{equation}
\lim_{r\rightarrow\infty}\frac{\log\sigma^{\prime}(-r)}{\log r}=-1.\label{eq:2-3}
\end{equation}

Step 5: Contradiction with the representation in Lemma \ref{l1}.
Suppose now that $\sigma$ is universally convex. Since $\sigma(0)=0$
and $\ensuremath{\sigma^{\prime}(0)=1},$ Lemma \ref{l1} then give
a probability measure $\mu$ on {[}0, 1{]} such that
\begin{equation}
\sigma^{\prime}(-r)=\exp\left[-2\int_{[0,1]}\log(1+rt)d\mu(t)\right].\label{eq:2-2}
\end{equation}
Taking logarithms, dividing by $\ensuremath{\log r},$ taking the
limit $r\rightarrow\infty$ by applying the dominated convergence
and employing \eqref{eq:2-3}, we have
\[
-1=-2\mu((0,1]),\quad\mu((0,1])=\frac{1}{2}.
\]
Indeed, $\log(1+rt)/\log r\rightarrow\mathbf{1}_{(0,1]}(t)$ as $r\rightarrow\infty$
and is uniformly bounded for $\ensuremath{r\geq2}.$ Since $\ensuremath{t\leq1},$
\eqref{eq:2-2} now implies
\[
\sigma^{\prime}(-r)\geq\exp[-2\mu((0,1])\log(1+r)]=\frac{1}{1+r}.
\]
Therefore, 
\[
\ensuremath{\liminf_{r\rightarrow\infty}Y(r)\geq1},
\]
 contradicting \eqref{eq:2-4}. Thus $\sigma$ is not universally
convex.

\section{\label{sec:A-p} A positivity lemma for $h$}

 In this section we prove a positivity lemma for the auxiliary function
$h$ defined in Subsection \ref{subsec:PT1}. Although Lemma \ref{l2}
is not needed for the proof of Theorem \ref{t1}, it is the key technical
ingredient in the proof of Theorem \ref{t2}: it implies the logarithmic
concavity condition $(\log g)''\leq0$ used in Subsection \ref{subsec:PT2}
to establish the universal starlikeness of the transformed function
$\widetilde{\sigma}(z).$
\begin{lem}
\label{l2} For all $x>0,$ we have 
\[
h''(x)+h'(x)<0.
\]
\end{lem}
\begin{proof}
It follows from \eqref{eq:3-5} that 
\[
2h(x)=\frac{U(x)}{(U(x)-1)^{3}}-\frac{u(x)}{(1-u(x))^{3}}.
\]
Differentiating both sides of the above identity yields
\[
2h'=-\frac{U(2U+1)}{(U-1)^{5}}+\frac{u(2u+1)}{(1-u)^{5}},
\]
 and 
\[
2h''=\frac{U(6U^{2}+8U+1)}{(U-1)^{7}}-\frac{u(6u^{2}+8u+1)}{(1-u)^{7}}.
\]
 Therefore 
\begin{equation}
2(h''+h')=\frac{U^{2}(-2U^{2}+9U+8)}{(U-1)^{7}}-\frac{u^{2}(-2u^{2}+9u+8)}{(1-u)^{7}}.\label{eq:3-2}
\end{equation}
 Put 
\[
s:=U-u>0,\qquad E:=e^{s}.
\]
Since $\ensuremath{U-\log U=u-\log u=x+1},$ we have $\log(U/u)=U-u=s,$
so $U/u=E.$ Thus, 
\begin{equation}
U=\frac{sE}{E-1},\qquad u=\frac{s}{E-1}.\label{eq:3-4}
\end{equation}
 Define 
\[
A:=(s-1)E+1,\qquad B:=E-1-s.
\]
 Then 
\[
A>0,\quad B>0
\]
and
\begin{equation}
U-1=\frac{A}{E-1},\qquad1-u=\frac{B}{E-1}.\label{eq:3-3}
\end{equation}
 Set 
\[
Q_{1}:=-2s^{2}+9s(E-1)+8(E-1)^{2},
\]
and 
\[
Q_{2}:=-2s^{2}E^{2}+9sE(E-1)+8(E-1)^{2}.
\]
 Substituting \eqref{eq:3-4} and \eqref{eq:3-3} into \eqref{eq:3-2}
yields 
\[
2(h''+h')=-\frac{s^{2}(E-1)^{3}}{A^{7}B^{7}}\,N(s),
\]
 where 
\[
N(s):=A^{7}Q_{1}-E^{2}B^{7}Q_{2}.
\]

It remains to prove $N(s)>0$ for all $\ensuremath{s>0}.$ We use
a finite exponential-polynomial positivity verification from \cite{BRS}.
Write 
\[
N(s)=f_{0}(s)=\sum_{j=0}^{11}P_{j}(s)e^{js}=\sum_{k\geq0}a_{k}s^{k},
\]
where $\{P_{j}(s)\}$ are polynomials in $s.$ For $k\geq0,$ let
$f_{k}$ denote the function obtained at stage $k$ of the reduction,
let $m_{k}$ be the largest exponential exponent in $f_{k},$ and
let $n_{k}$ be the degree of the constant-exponent term in $f_{k}.$
Define
\[
f_{k+1}:=e^{-s}f_{k}^{(n_{k}+1)}.
\]
 Then $m_{k+1}=m_{k}-1$ and the degrees $n_{k}$ are exactly the
numbers listed in the second column of Table \ref{tab:positivity-certificate}.
At stage $k$ we check 
\begin{equation}
f_{k}^{(r)}(0)\geq0,\quad r=0,1,\ldots,n_{k}.\label{eq:2-7}
\end{equation}
Thus there are $n_{k}+1$ checks at stage $k.$ Table \ref{tab:positivity-certificate}
records these checks, separated into zero checks and positive checks;
no negative check occurs. The last row gives the totals 23 zero checks
and 63 positive checks.

\begin{table}[htbp] \centering \renewcommand{\arraystretch}{1.15} \begin{tabular}{c|c|c|c} \hline Stage \(k\) & Degree \(n_k\) & Zero checks & Positive checks \\ \hline 0  & 2 & 3 & 0  \\ 1  & 3 & 4 & 0  \\ 2  & 7 & 8 & 0  \\ 3  & 8 & 8 & 1  \\ 4  & 9 & 0 & 10 \\ 5  & 8 & 0 & 9  \\ 6  & 8 & 0 & 9  \\ 7  & 9 & 0 & 10 \\ 8  & 8 & 0 & 9  \\ 9  & 7 & 0 & 8  \\ 10 & 3 & 0 & 4  \\ 11 & 2 & 0 & 3  \\ \hline Total & & 23 & 63 \\ \hline \end{tabular} \caption{Exact finite integer checks for the positivity certificate of \(N(s)\).} \label{tab:positivity-certificate} \end{table}The
reduction principle is used as follows. Suppose that all Taylor coefficients
of $f_{k+1}$ are nonnegative. Since 
\[
f_{k}^{(n_{k}+1)}(s)=e^{s}f_{k+1}(s)
\]
and $e^{s}$ has positive Taylor coefficients, all Taylor coefficients
of $f_{k}^{(n_{k}+1)}$ are nonnegative. This implies that the Taylor
coefficients of $f_{k}$ of order at least $n_{k}+1$ are nonnegative.
If, in addition, the checks \eqref{eq:2-7} hold, then the first $n_{k}+1$
Taylor coefficients of $f_{k}$ are also nonnegative. Hence all Taylor
coefficients of $f_{k}$ are nonnegative. Iterating this implication
backwards from the final polynomial to $f_{0}=N$ we conclude that
all Taylor coefficients of $N$ are nonnegative, provided that the
final polynomial has nonnegative coefficients. The final polynomial
is checked separately in lines 14\textendash 16 of the Maple verifier
in Appendix \ref{app:verifier}; its coefficients are nonnegative
and at least one of them is positive. Actually, in the reduction chain
one has the exact identity
\[
f_{3}^{(8)}(0)=N^{(23)}(0).
\]
This is the unique positive check in the row $k=3$ of Table \ref{tab:positivity-certificate},
namely the check with $r=8$ in \eqref{eq:2-7}. The standalone Maple
verifier gives
\[
\ensuremath{N^{(23)}(0)=14841922455360000>0}.
\]
Since 
\[
N^{(23)}(0)=23!a_{23},
\]
where $a_{23}$ is the coefficient of $s^{23}$ in $N(s),$ this positive
value proves that $a_{23}>0.$ Thus, at least one Taylor coefficient
of $N$ is positive. Together with the nonnegativity of all Taylor
coefficients of $N,$ this yields
\[
N(s)>0\qquad(s>0).
\]
Therefore,
\[
h''(x)+h'(x)<0\qquad(x>0).
\]
This proves Lemma \ref{l2}.
\end{proof}

\section{\label{sec:Univ} Universal convexity of the generalized Ramanujan
generating function}

\subsection{Definition and main result for $\sigma_{\alpha}$}

Let
\[
\sigma_{\alpha}(z):=\frac{1}{\theta(1)-1/3}\sum_{n\geq1}(\theta(n)-1/3)\frac{z^{n}}{n^{\alpha}}.
\]
In \cite[Theorem 3]{BRS}, Bakan, Ruscheweyh and Salinas proved that
$\sigma_{\alpha}(z)$ is universally starlike for every $\alpha\geq1,$
which seems to be the first result placing the Ramanujan sequence
into the study of analytic univalent functions. This is also an important
step towards the conjecture of Ruscheweyh, Salinas and Sugawa.

It is a natural question whether $\sigma_{\alpha}(z)$ is universally
convex for every $\alpha\geq1.$ In this paper, we provide a sharp
criterion of universal convexity for $\sigma_{\alpha}(z),$ answering
this question positively and giving a general extension of Theorem
\ref{t1}.
\begin{thm}
\label{t2} The function $\sigma_{\alpha}(z)$ is universally convex
if and only if $\alpha\geq1.$
\end{thm}
Taking $\alpha=0$ recovers the non-universal-convexity stated in
Theorem \ref{t1}, while the case $\alpha\geq1$ gives a positive
and sharp result.

\subsection{\label{subsec:PT2} Proof of Theorem \ref{t2}}

We first prove the sufficiency.  We split the proof into two steps.

Step 1: Universal starlikeness of $\sigma.$ It follows from \eqref{eq:2-5}
that 
\[
\frac{\sigma(z)}{z}=\frac{1}{c}\int_{0}^{\infty}\frac{h(x)}{e^{x}-z}\,dx=\frac{1}{c}\int_{0}^{1}\frac{h(-\log t)}{1-tz}\,dt.
\]
 Consider the involution 
\[
\widetilde{\sigma}(z)=-\sigma\!\left(\frac{z}{z-1}\right).
\]
It is analytic and normalized on $\ensuremath{\Lambda}$ since $\widetilde{\sigma}(0)=0,\;\widetilde{\sigma}'(0)=1.$
Then
\[
\frac{\widetilde{\sigma}(z)}{z}=\frac{1}{c}\int_{0}^{1}\frac{h(-\log(1-t))}{1-tz}\,dt.
\]
 Take $t\mapsto1/t.$ Then 
\begin{equation}
\frac{\widetilde{\sigma}(z)}{z}=\int_{1}^{\infty}\frac{g(\log t)/t}{t-z}\,dt,\label{eq:4-1}
\end{equation}
where 
\[
g(x)=\frac{h(\kappa(x))}{c},\quad\kappa(x)=-\log(1-e^{-x}).
\]
 We want to apply \cite[Lemma 2]{BRS}. It suffices to check 
\[
g(0+)=0,\qquad\int_{0}^{\infty}e^{-x}g(x)\,dx=1,\qquad g'\ge0,\qquad(\log g)''\le0.
\]
 The first three conditions follow from $\ensuremath{h(\infty)=0},$
$\ensuremath{h'<0},$ and $\int_{0}^{\infty}e^{-x}h(x)\,dx=c.$ For
the final condition, put 
\[
b(y):=\log h(y),\qquad q:=e^{\kappa(x)}-1>0.
\]
 Since 
\[
\kappa'(x)=-q,\qquad\kappa''(x)=q(q+1),
\]
we obtain 
\[
(\log g)''(x)=q^{2}\bigl(b''(\kappa(x))+b'(\kappa(x))\bigr)+q\,b'(\kappa(x))<0.
\]
Indeed, $\ensuremath{b'=h'/h<0},$ and Lemma \ref{l2} gives 
\[
b''+b'=\frac{h''+h'}{h}-\left(\frac{h'}{h}\right)^{2}<0.
\]
Thus, by \cite[Lemma 2]{BRS}, there exists a probability measure
$\mu$ on {[}0,1{]} such that
\[
\frac{\widetilde{\sigma}(z)}{z}=\exp L_{\mu}(z).
\]
Together with \cite[Theorem 5]{BRS}, this shows that $\widetilde{\sigma}$
is universally starlike.

Define
\[
Tf(z):=-f\!\left(\frac{z}{z-1}\right).
\]
 Then Theorems 1.4 and 1.5 in \cite{RSS} imply that if $f$ is universally
starlike, so is $Tf(z).$ Since $\ensuremath{T^{2}f=f},$ the universal
starlikeness of $\widetilde{\sigma}$ implies that of $\ensuremath{\sigma}=\sigma_{0}.$

Step 2: Universal convexity of $\sigma_{\alpha}$ for $\ensuremath{\alpha\ge1}.$
Now fix $\alpha\ge1$ and put $\ensuremath{\beta:=\alpha-1\ge0}.$
The universal starlikeness of $\sigma$ shows that $\sigma_{\beta}$
is universally starlike. Hence, by \cite[Theorem 5]{BRS},
\[
\sigma_{\alpha}'(z)=\frac{\sigma_{\beta}(z)}{z}=\exp L_{\mu_{\beta}}(z)
\]
 for some probability measure $\mu_{\beta}$ on $\ensuremath{[0,1]}.$
Let 
\[
\nu_{\alpha}=\frac{\mu_{\beta}+\delta_{0}}{2}.
\]
Since $\ensuremath{L_{\delta_{0}}=0},$ we have $2L_{\nu_{\alpha}}=L_{\mu_{\beta}}+L_{\delta_{0}}=L_{\mu_{\beta}}.$
Then
\[
\sigma_{\alpha}'(z)=\exp\bigl(2L_{\nu_{\alpha}}(z)\bigr).
\]
 Lemma \ref{l1} therefore proves universal convexity of $\sigma_{\alpha}$
for every $\ensuremath{\alpha\ge1},$ including the endpoint $\ensuremath{\alpha=1}.$

We next prove the necessity.  We divide the proof into two steps.

(1) Exclusion of $0\le\alpha<1.$ Define the fractional integral 
\[
v_{\alpha}(x):=\frac{1}{\Gamma(\alpha)}\int_{0}^{x}(x-y)^{\alpha-1}h(y)\,dy,\qquad0<\alpha<1
\]
and $v_{0}(x):=h(x).$ Tonelli's theorem and the gamma integral yield,
for $\ensuremath{n\ge1},$ 
\[
\int_{0}^{\infty}e^{-nx}v_{\alpha}(x)\,dx=\frac{a_{n}}{n^{\alpha}}.
\]
 Consequently, for $\ensuremath{0\le\alpha<1},$ 
\[
\sigma_{\alpha}(z)=\frac{1}{c}\int_{0}^{\infty}\frac{zv_{\alpha}(x)}{e^{x}-z}\,dx,\qquad\sigma_{\alpha}'(-r)=\frac{1}{c}\int_{0}^{\infty}\frac{e^{x}v_{\alpha}(x)}{(e^{x}+r)^{2}}\,dx>0.
\]
 The first equality holds initially for $|z|<1$ and extends analytically
to $\ensuremath{\Lambda};$ its derivative gives the second equality.

For each fixed $\ensuremath{0<\alpha<1},$ there are positive constants
$C_{\alpha},D_{\alpha}$ such that, for all sufficiently large $\ensuremath{x},$
\[
C_{\alpha}x^{\alpha-1}\le v_{\alpha}(x)\le D_{\alpha}x^{\alpha-1}.
\]
For the lower bound, use $\ensuremath{(x-y)^{\alpha-1}\ge x^{\alpha-1}}.$
For the upper bound, split it at $\ensuremath{y=x/2}.$ On $\ensuremath{[0,x/2]},$
use integrability of $h$ and $\ensuremath{(x-y)^{\alpha-1}\le(x/2)^{\alpha-1}}.$
On $\ensuremath{[x/2,x]},$ use $h(y)\le Cx^{-2}$ from \eqref{eq:2-6},
giving an upper bound of order $\ensuremath{x^{\alpha-2}}.$ Near
$\ensuremath{0},$ boundedness of $h$ gives $\ensuremath{v_{\alpha}(x)\le\|h\|_{\infty}x^{\alpha}/\Gamma(\alpha+1)}.$
Thus every $v_{\alpha}$ with $0\le\alpha<1$ is bounded and tends
to $0$ at infinity. Moreover, for sufficiently large $\ensuremath{x},$
\begin{equation}
v_{\alpha}(x)\ge Cx^{-q},\qquad q=\begin{cases}
2, & \mathrm{if}\;\alpha=0,\\
1-\alpha, & \mathrm{if}\;0<\alpha<1,
\end{cases}\label{eq:4-2}
\end{equation}
 where $C>0$ may depend on $\ensuremath{\alpha}.$

Put 
\[
L=\log r,\qquad s=x-L,\qquad x=L+s=\log r+s.
\]
 Then $dx=ds$ and 
\[
e^{x}=e^{L+s}=re^{s},
\]
 so 
\[
(e^{x}+r)^{2}=(re^{s}+r)^{2}=r^{2}(e^{s}+1)^{2}.
\]
 Hence 
\[
\frac{e^{x}}{(e^{x}+r)^{2}}\,dx=\frac{re^{s}}{r^{2}(e^{s}+1)^{2}}\,ds=\frac{1}{r}\cdot\frac{e^{s}}{(1+e^{s})^{2}}\,ds=\frac{1}{r}K(s)\,ds,
\]
 where 
\[
K(s)=\frac{e^{s}}{(1+e^{s})^{2}}.
\]
 Therefore 
\[
\sigma_{\alpha}'(-r)=\frac{1}{c}\int_{-L}^{\infty}K(s)v_{\alpha}(L+s)\,\frac{ds}{r}.
\]
 Namely, 
\[
r\sigma_{\alpha}'(-r)=\frac{1}{c}\int_{-L}^{\infty}K(s)v_{\alpha}(L+s)\,ds.
\]
 Since $v_{\alpha}$ is extended by $0$ on the negative half-line,
the integrand vanishes for $\ensuremath{s<-L}.$ Thus 
\[
r\sigma_{\alpha}'(-r)=\frac{1}{c}\int_{\mathbb{R}}K(s)v_{\alpha}(L+s)\,ds.
\]
 Since $K$ is integrable, $\ensuremath{\int_{\mathbb{R}}K=1},$ and
$v_{\alpha}$ is bounded with limit $0$ at infinity, the dominated
convergence shows that
\begin{equation}
r\sigma_{\alpha}'(-r)\to0\qquad(r\to\infty).\label{eq:4-4}
\end{equation}
 On the other hand, restricting to $0\le s\le1$ and using \eqref{eq:4-2}
gives 
\[
C_{1}(\log r)^{-q}\le r\sigma_{\alpha}'(-r)\le C_{2}
\]
 for all sufficiently large $\ensuremath{r},$ with $\ensuremath{C_{1},C_{2}>0}.$
Taking logarithms and dividing by $\log r$ therefore gives 
\begin{equation}
\lim_{r\to\infty}\frac{\log\sigma_{\alpha}'(-r)}{\log r}=-1.\label{eq:4-3}
\end{equation}

Suppose $\sigma_{\alpha}$ is universally convex for some $\ensuremath{0\le\alpha<1}.$
Since the $\alpha=0$ case is included in Theorem \ref{t1} and was
treated in Subsection \ref{subsec:PT1}, we assume that $0<\alpha<1.$
By Lemma \ref{l1}, there exists a probability measure $\mu$ on $[0,1]$
such that
\[
\log\sigma_{\alpha}'(-r)=-2\int_{[0,1]}\log(1+rt)\,d\mu(t).
\]
 For $\ensuremath{r\ge2},$ the quotient $\log(1+rt)/\log r$ is bounded
uniformly in $t\in[0,1]$ and converges to $1$ for $t>0$ and to
$0$ for $\ensuremath{t=0}.$ The dominated convergence and \eqref{eq:4-3}
imply 
\[
-1=-2\mu((0,1]),\qquad\text{so}\qquad\mu((0,1])=\frac{1}{2}.
\]
 Since $\log(1+rt)\le\log(1+r)$ on $\ensuremath{[0,1]},$ we get
\[
\sigma_{\alpha}'(-r)\geq\exp\left(-2\int_{(0,1]}\log(1+r)\,d\mu(t)\right)=\frac{1}{1+r}.
\]
 Thus, 
\[
\liminf_{r\to\infty}r\sigma_{\alpha}'(-r)\ge1,
\]
 contradicting \eqref{eq:4-4}. Thus $\sigma_{\alpha}$ is not universally
convex for any $\ensuremath{0\le\alpha<1}.$

(2) Exclusion of $\alpha<0.$ Suppose $\alpha<0$ and $\sigma_{\alpha}$
is universally convex. Choose 
\[
m=\lceil-\alpha\rceil\ge1,
\]
so that $\ensuremath{0\le\alpha+m<1}.$ By \cite[Theorem 2]{BRS},
dividing its $n$th coefficient by $n^{m}$ makes $\sigma_{\alpha+m}$
universally convex, contradicting the fact that $\sigma_{\beta}$
is not universally convex for any $\ensuremath{0\le\beta<1}.$ This
excludes every negative $\alpha.$ Thus, $\alpha\geq1.$ This completes
the proof of Theorem \ref{t2}. \qed

\appendix
\section{Exact Maple verifier for Lemma \ref{l2}} \label{app:verifier}
The following standalone Maple script provides an exact verification of the finite positivity certificate used in the proof of Lemma \ref{l2}. It uses exact integer and rational arithmetic only, and reproduces the checks recorded in Table \ref{tab:positivity-certificate}. Running the script in Maple should end with
\begin{verbatim} Exact certificate verified: N^(23)(0)=14841922455360000>0 \end{verbatim}
\noindent The complete script is reproduced below.
\begin{lstlisting}[language=Maple,   caption={Exact verification for Lemma \ref{l2}},   label={lst:maple-verifier}] restart:
# ============================================================ # Exact verifier for the positivity certificate in Lemma 3.1 # ============================================================
Dpow := proc(P, j, m) local r, k;   r := expand(P);   for k to m do     r := expand(diff(r, s) + j*r);   end do;   r; end proc:
ReduceExpPoly := proc(Plist, M)   local n, Q, e;   n := degree(expand(Plist[1]), s);   Q := [seq(0, e=1..M)];   for e from 1 to M do     Q[e] := expand(Dpow(Plist[e+1], e, n+1));   end do;   Q; end proc:
E := 'E': A := (s-1)*E + 1: B := E - 1 - s: Q1 := -2*s^2 + 9*s*(E-1) + 8*(E-1)^2: Q2 := -2*s^2*E^2 + 9*s*E*(E-1) + 8*(E-1)^2: Npoly := expand(A^7*Q1 - E^2*B^7*Q2):
M := degree(Npoly, E); Plist := [seq(expand(coeff(Npoly, E, j)), j=0..M)];
f := Plist: m := M: k := 0:
while m > 0 do   n := degree(expand(f[1]), s);   fexpr := add(f[j+1]*exp(j*s), j=0..m);
  zeroChecks := 0:   posChecks := 0:   negChecks := 0:
  for r from 0 to n do     if r = 0 then       val := expand(eval(fexpr, s=0));     else       val := expand(eval(diff(fexpr, s$r), s=0));     end if;
    if val = 0 then       zeroChecks := zeroChecks + 1;     elif val > 0 then       posChecks := posChecks + 1;     else       negChecks := negChecks + 1;     end if;   end do;
  printf("Stage %d: n_k = %d, zero checks = %d, positive checks = %d, negative checks = %d\n",          k, n, zeroChecks, posChecks, negChecks);
  if negChecks > 0 then     error "Negative derivative check at stage %d", k;   end if;
  f := ReduceExpPoly(f, m);   m := m - 1;   k := k + 1; end do:
finalPoly := expand(f[1]); degFinal := degree(finalPoly, s); coeffList := [seq(coeff(finalPoly, s, r), r=0..degFinal)];
if min(coeffList[]) < 0 then   error "Final polynomial has a negative Taylor coefficient"; end if:
if max(coeffList[]) <= 0 then   error "Final polynomial has no positive Taylor coefficient"; end if:
printf("Final polynomial: %a\n", finalPoly); printf("All Taylor coefficients of the final polynomial are nonnegative, and at least one is positive.\n");
N := subs(E=exp(s), Npoly): cert := expand(eval(diff(N, s$23), s=0)):
printf("N^(23)(0) = %a\n", cert);
if cert = 14841922455360000 then   print("Exact certificate verified: N^(23)(0) = 14841922455360000 > 0"); else   error "Certificate mismatch: got %a", cert; end if; \end{lstlisting}

\section*{Acknowledgement}

 This work was partially supported by the Scientific Research Fund
of Hunan Provincial Education Department (Grant No. 25B0010).


\begin{thebibliography}{1}
\bibitem{AJ}J.A. Adell and P. Jodrá, On the complete monotonicity
of a Ramanujan sequence connected with $e^{n},$ Ramanujan J. 16 (2008),
1\textendash 5.

\bibitem{BRS}A. Bakan, S. Ruscheweyh and L. Salinas, On geometric
properties of the generating function for the Ramanujan sequence.
Ramanujan J. 46(1)(2018), 173\textendash 188.

\bibitem{FGKP}P. Flajolet, P.J. Grabner, P. Kirschenhofer and H.
Prodinger, On Ramanujan's Q-function. J. Comput. Appl. Math. 58(1)(1995),
103\textendash 116.

\bibitem{K}S. Koumandos, A Bernstein function related to Ramanujan's
approximations of $exp(n).$ Ramanujan J. 30(3)(2013), 447\textendash 459.

\bibitem{R}S. Ramanujan, Collected Papers. Chelsea, New York, 1962.

\bibitem{RSS}S. Ruscheweyh, L. Salinas and T. Sugawa, Completely
monotone sequences and universally prestarlike functions. Israel J.
Math. 171 (2009), 285\textendash 304.

\bibitem{S}G. Szegö, Über einige von S. Ramanujan gestellte Aufgaben.
J. Lond. Math. Soc. 3(1928), 225\textendash 232.

\bibitem{W}G.N. Watson, Theorems stated by Ramanujan (V): approximations
connected with $e^{x}.$ Proc. Lond.Math. Soc. 29(1929), 293\textendash 308.
\end{thebibliography}
\end{document}